\documentclass[12pt]{article}

\usepackage{a4wide}
\usepackage{amsmath,amssymb,amsthm}
\usepackage{color,graphicx}

\newtheorem{theorem}{Theorem}

\newtheorem{proposition}{Proposition} 

\theoremstyle{definition}
\newtheorem*{example}{Example}
\newtheorem{remark}{Remark}

\newenvironment{dense_itemize}{%
  \begin{list}{$\triangleright$}%
    {\setlength{\topsep}{1mm}%
      \setlength{\partopsep}{0mm}%
      \setlength{\parskip}{0mm}%
      \setlength{\parsep}{0mm}%
      \setlength{\itemsep}{0mm}%
      \setlength{\labelwidth}{4mm}%
      \setlength{\leftmargin}{0mm}%
      \addtolength{\leftmargin}{\labelwidth}%
      \addtolength{\leftmargin}{\labelsep}%
      \setlength{\itemindent}{0mm}}}%
  {\end{list}}

\newcommand{\F}{{\mathbb F}}
\newcommand{\R}{{\mathbb R}}

\newcommand{\C}{{\mathcal C}}
\renewcommand{\S}{{\mathcal S}}
\newcommand{\T}{{\mathcal T}}
\renewcommand{\phi}{\varphi}

\DeclareMathOperator{\st}{star}

\begin{document}

\title{Many Triangulated 3-Spheres}
\author{\textsc{Julian Pfeifle}\thanks{Supported by the Deutsche Forschungsgemeinschaft
    within the European graduate program \emph{Combinatorics, Geometry, and
      Computation} (No. GRK 588/1) and an MSRI post-doctoral fellowship}
  \quad and\quad\setcounter{footnote}{6}%
  \textsc{G\"unter M.~Ziegler}\thanks{Partially supported by Deutsche
  Forschungs-Gemeinschaft (DFG), via the
  DFG Research Center ``Mathematics in the Key Technologies'' (FZT86),
  the Research Group ``Algorithms, Structure, Randomness'' (Project ZI 475/3),
  and a Leibniz grant (ZI 475/4)}\\
  TU Berlin, MA 6-2\\
  \texttt{$\{$pfeifle,ziegler$\}$@math.tu-berlin.de}}
\date{}

\maketitle

\begin{abstract}
  We construct $2^{\Omega(n^{5/4})}$ combinatorial types of
  triangulated $3$-spheres on $n$~vertices. Since by a result of
  Goodman and Pollack (1986) there are no more than $2^{O(n\log n)}$
  combinatorial types of simplicial $4$-polytopes, this proves that
  asymptotically, there are far more combinatorial types of
  triangulated $3$-spheres than of simplicial $4$-polytopes on
  $n$~vertices. This complements results of Kalai (1988), who had
  proved a similar statement about $d$-spheres and $(d+1)$-polytopes
  for fixed~$d\ge4$.


\end{abstract}

\section{Introduction}

In 1988, Kalai~\cite{Kalai88} proved a lower bound of
\[ 
    \log s(d,n) \ = \ \Omega(n^{\lfloor d/2\rfloor}) \qquad 
    \text{for fixed } d\ge3
\]
for the number $s(d,n)$ of distinct combinatorial types of simplicial
PL $d$-spheres on $n$~vertices.\footnote{Here and in the following, 
we use Landau's $O$-notation for positive functions:
$f = \Omega(g)$ denotes that there is a positive constant
$c$ such that $f(n)\ge c g(n) $ holds for all sufficiently large $n$.
Similarly $f = O(g)$ if there is a $c'>0$ such that
$f(n)\le c' g(n)$ for all large $n$,
and $f= \Theta(g)$ denotes that both conditions hold.}
Combining this with
Goodman and Pollack's~\cite{Goodman-Pollack86,Goodman-Pollack87} upper
bound
\[
    \log p(d,n) \;\le\; d(d+1)n\log n
\]
for the number $p(d,n)$ of combinatorial types of simplicial
$d$-polytopes on $n$ vertices, he derived that for $d\ge4$, there are
far more simplicial $d$-spheres than simplicial $(d+1)$-polytopes. In
particular, \emph{most} of these spheres, in a very strong sense, are
not polytopal, i.e.~there is no convex polytope with the same face
lattice.  On the other hand, we proved in earlier
work~\cite{Pfeifle02} that in dimension $d=3$, Kalai's construction
produces only polytopal spheres, and up to now only few families of
non-polytopal $3$-spheres were known.

In this paper, 
we combine two constructions from a recent paper
by Eppstein, Kuperberg \& Ziegler \cite{EKZ02} to 
show for the first time that for $n$ large enough,
there are far more simplicial $3$-spheres than $4$-polytopes on
$n$~vertices. 

\begin{theorem}\label{thm:many-spheres}
  There are at least 
  \[
      s(3,n) \ = \ 2^{\displaystyle \Omega(n^{5/4})}
  \]
  combinatorially non-isomorphic simplicial $3$-spheres on $n$
  vertices.
\end{theorem}

In brief, we prove Theorem~\ref{thm:many-spheres} by producing a
cellular decomposition~$\S$ of~$S^3$ with $n$~vertices and
$\Theta(n^{5/4})$ octahedral facets, and triangulating each octahedron
independently.  The cellulation~$\S$ is constructed from a Heegaard
splitting $S^3=H_1\cup H_2$ of~$S^3$ of high genus by appropriately
subdividing the thickened boundary surface $(H_1\cap H_2)\times[0,1]$.

Because of their sheer number, \emph{most} of the spheres we
construct are combinatorially distinct: There can be at most
$n!$~spheres combinatorially isomorphic to any given one, where
$n!<n^n=2^{n\log n}$. Also note that the only currently known upper bound
for~$s(3,n)$ is the rather crude estimate $\log s(3,n)=O(n^2\log n)$
obtained from Stanley's proof of the Upper Bound Theorem for
spheres~\cite{Stanley75}.

\section{Background}

For $d\le2$ all simplicial $d$-spheres are realizable as polytopes:
$1$-dimensional spheres are trivial to realize, and Steinitz' famous
theorem~\cite{Steinitz1922},~\cite{Steinitz-Rademacher34} from the
beginning of the 20th century asserts that \emph{all} $2$-spheres,
including the non-simplicial ones, are polytopal (i.e., they arise as
boundary complexes of $3$-dimensional polytopes). Tutte~\cite{Tutte80}
showed in 1980 that the number of combinatorially distinct rooted
simplicial $3$-polytopes with $n$~vertices is asymptotically
\[
    \frac{3}{16\sqrt{6\pi n^5}}\,\left(\frac{256}{27}\right)^{n-2},
\]
and Bender~\cite{Bender87} established sharp asymptotic formulas counting the
number of unrooted $3$-dimensional polytopes.


The first example---the so-called \emph{Br\"uckner sphere}---of a
simplicial sphere that is \emph{not} the boundary
complex of a polytope was inadvertedly found by
Br\"uckner~\cite{Brueckner1910} in 1910 in an attempt to enumerate all
combinatorial types of $4$-polytopes with $8$~facets. As noted in 1967
by Gr\"unbaum and Sreedharan~\cite{Gruenbaum-Sreedharan67}, one of the
$3$-dimensional
complexes that Br\"uckner thought to represent a polytope is in fact
\emph{not} realizable in a convex way in~$\R^4$. As the (polytopal)
complex Br\"uckner considered is simple (any vertex is contained in
exactly $4$~facets), its combinatorial dual is a simplicial
$3$-sphere.  

Another known `sporadic' example of a non-polytopal simplicial
sphere is \emph{Barnette's sphere}~\cite{Barnette70}, which is nicely
explained in~\cite[Chapter~III.4]{Ewald96}. From these two examples
one can build infinite series, but apart from such sporadic families,
no substantial number of non-polytopal spheres on a fixed number of
vertices was known until Kalai's 1988 construction.

The related problem of estimating the number $t(3,m)$ of combinatorial
types of simplicial $3$-spheres with $m$~\emph{facets} has attracted
attention in gravitational quantum physics~\cite{DJ95}. 
Gromov~\cite{Gromov2000} has 
asked whether there exists a constant~$c>0$ such that $t(3,m)\le 2^{cm}$. 
By duality, this is equivalent to bounding the number of
\emph{simple} $3$-spheres on $m$~vertices.
The problem is ``dual'' to the one we treat here,
but seems to require different methods.

\section{Definitions and notation}

A \emph{cellulation}~$\C$ of a manifold $X$ is a finite CW complex
whose underlying space is~$X$. $\C$~is \emph{regular} if all closed
cells are embedded, and \emph{strongly regular} if in addition the
intersection of any two cells is a cell. The \emph{star} of a
cell~$\sigma\in\C$ is the union of the closure of all cells
containing~$\sigma$, and the \emph{link} of~$\sigma$ consists of all
cells of $\st\sigma$ not incident to~$\sigma$. The entry~$f_i$ of the
\emph{$f$-vector} $f(\C)=(f_0,f_1,\dots)$ of a cellulation counts the
number of $i$-dimensional cells. The $d$-dimensional cells are called
\emph{facets}, and $(d-1)$-dimensional ones \emph{ridges}.

\section{The ingredients for the construction}

\subsection{Heffter's embedding of the complete graph}

In 1898, Heffter~\cite{Heffter1898} constructed remarkable
cellulations of closed orientable surfaces:

\begin{proposition}\label{prop:heffter}
  Let $q=4k+1$ be a prime power, and $\alpha$ be any generator of the
  cyclic group~$\F_q^*$ of invertible elements of the finite
  field~$\F_q$ on $q$~elements. Then there exists a regular but not
  strongly regular cellulation~$\C_q^\alpha$ of the closed orientable
  surface~$S_g$ of genus $g=q(q-5)/4+1$ with $f$-vector
  $\big(q,{q\choose 2},q\big)$, all of whose $2$-cells are
  $(q-1)$-gons. $\C_q^\alpha$~can be refined to a strongly regular
  triangulation~$\T_q^\alpha$ of~$S_g$ with $f$-vector
  $\big(2q,{q\choose 2}+q(q-1),q(q-1)\big)$.
\end{proposition}

\begin{proof}
There exist infinitely many \emph{prime}
numbers~$q$ of the form $q=4k+1$; see~\cite{Erdoes35}. For any prime
\emph{power}~$q$ of this form, take as vertices of the
cellulation~$\C^\alpha_q$ the elements of~$\F_q$, and as $2$-cells the
$(q-1)$-gons (compare Figure~\ref{fig:Heffter-2}, left)
\[
    F^\alpha(s) \ =\  
    \Big\{ v^\alpha(s,k) =
    s+\frac{\alpha^k-1}{\alpha-1}:0\le k\le q-2\Big\} \qquad
    \text{for } s\in\F_q.
\]

\begin{figure}[htbp]
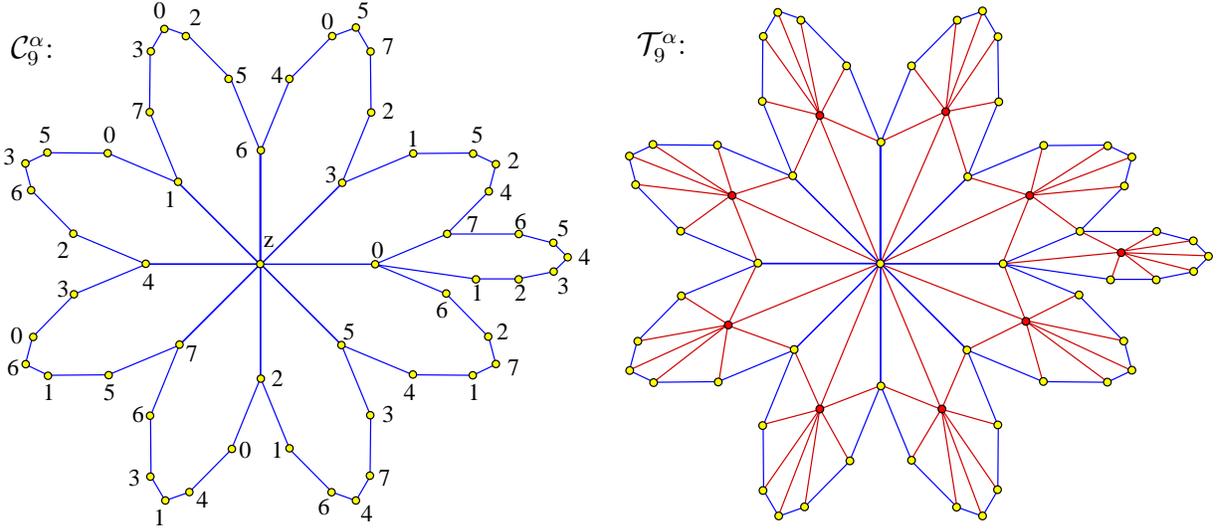
 
  \centering
  \begin{minipage}[c]{.45\linewidth}
    \input{heffter-2.pstex_t}
  \end{minipage} 
  \qquad
  \begin{minipage}[c]{.45\linewidth}
    \input{heffter-2-subdiv.pstex_t}
  \end{minipage}
\caption{\emph{Left:} The Heffter cellulation~$\C_9^\alpha$ of a
  surface~$S_g$ of genus $g=10$ for $\alpha=2x+2\in\F_9\protect\cong
  \F_3[x]/ \langle x^2+x+2\rangle$.  The vertex $z$ corresponds to
  $0\in\F_9$, and the vertices labeled $i$ to the element $\alpha^i$.
  Note that any two of the $q=9$ vertices are adjacent, and that all
  vertices in any given one of the $9$~polygons are distinct.
  However, the link of every vertex contains identified vertices, and
  so the vertex stars are not embedded.  \emph{Right:} After
  subdividing the surface to the triangulation~$\T_q^\alpha$ using $q$
  new vertices, all stars are embedded disks}
  \label{fig:Heffter-2}  
\end{figure} 

It is straightforward to check (see~\cite{Heffter1898} and~\cite[Lemma
12]{EKZ02}) that this cellulation is regular (all vertices in
each~$F(s)$ are distinct), neighborly (any two vertices are connected
by an edge), and closed (any edge is shared by exactly two polygons),
but not strongly regular (any two polygons share $q-2$ vertices). An
Euler characteristic calculation yields the genus of the underlying
surface~$S_g$ of~$\C_q^\alpha$.  By subdividing each polygon as in
Figure~\ref{fig:Heffter-2} (right), the cellulation becomes strongly
regular with the stated $f$-vector.  
\end{proof}

\begin{remark}
  This cellulation was independently obtained in~\cite{EKZ02} as an
  abelian covering of the canonical one-vertex cellulation of~$S_g$.
\end{remark}

\begin{remark}
  Heffter's original construction involved only prime numbers. As it
  turns out, allowing prime powers becomes necessary for symmetric
  embeddings: According to Biggs~\cite{Biggs69}, if the complete
  graph~$K_n$ embeds into a closed orientable surface in a symmetric
  way (i.e. there exists a ``rotary'' or ``chiral'' combinatorial
  automorphism, see~\cite{Wilson2002}), then $n$~is the power of a
  prime number, and James \& Jones~\cite{JamesJones85} showed that
  \emph{any} such embedding of~$K_n$ is actually one from Heffter's
  family.
\end{remark}

\begin{remark}
  Two cellulations $\C^\alpha_q$ and $\C^\beta_q$ are combinatorially
  distinct for $\beta\ne\alpha,1/\alpha\in\F_q$: By~\cite{Heffter1898}, the
  only automorphisms of~$\C^\alpha_q$ are induced by affine maps
  $\phi:\F_q\to\F_q$, $x\mapsto ax+b$ with $a\in\F_q^*,b\in\F_q$. An
  easy calculation shows that requiring
  $\phi(v^\alpha(s,k+i))=v^\beta(t,\ell+i)$
  resp.~$\phi(v^\alpha(s,k+i))=v^\beta(t,\ell-i)$ for $t\in\F_q$,
  $0\le\ell\le q-2$ and $i=0,1,2$ already implies $\beta=\alpha$
  resp.~$\beta=1/\alpha$. 
\end{remark}

\subsection{The E-construction}

\begin{proposition} \label{prop:e-construction} \textup{\cite{EKZ02}}
  Given a cellulation~$\C$ of a $d$-dimensional manifold~$M$ with
  boundary with $f$-vector $(f_0,f_1,\dots,f_d)$ and $f_{d-1}^{\rm
    in}$ interior ridges, there exists a cellulation~$E(\C)$ of~$M$
  with $f_0+f_d$ vertices consisting of $f_{d-1}^{\rm in}$ bipyramids
  and $f_{d-1}-f_{d-1}^{\rm in}$ pyramids.
\end{proposition}
 
\begin{proof}
Cone a new vertex to the inside of each
$d$-cell $F$ of~$M$ to create $f_{d-1}+f_{d-1}^{\rm in}$ pyramids,
then combine each pair of pyramids over the same interior ridge to a
bipyramid.  
\end{proof}

\begin{figure}[htbp]
  \centering
  \parbox[c]{.4\linewidth}{\includegraphics[width=\linewidth]{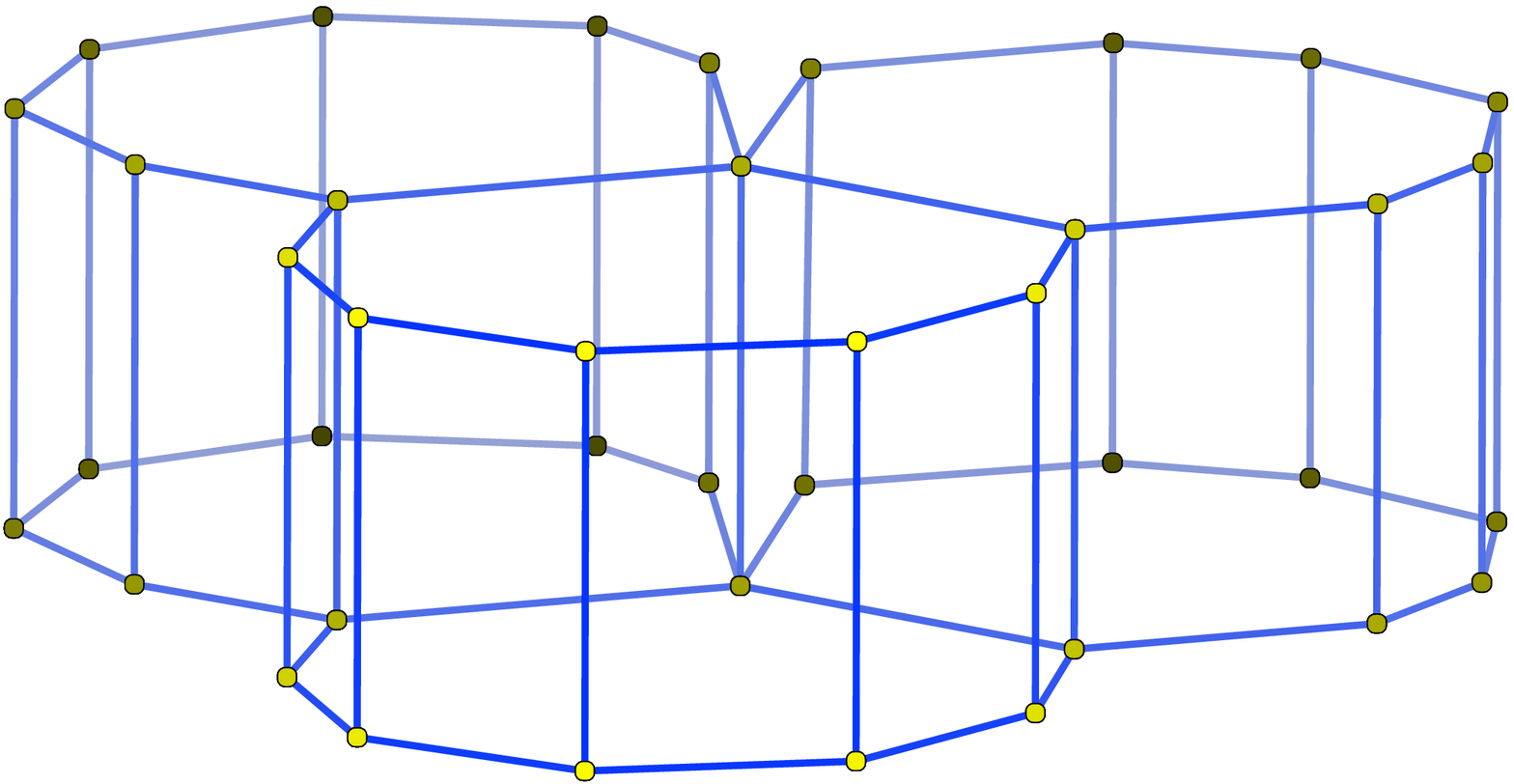}}
  \parbox[c]{.1\linewidth}{\quad$\longrightarrow$\quad}
  \parbox[c]{.4\linewidth}{\includegraphics[width=\linewidth]{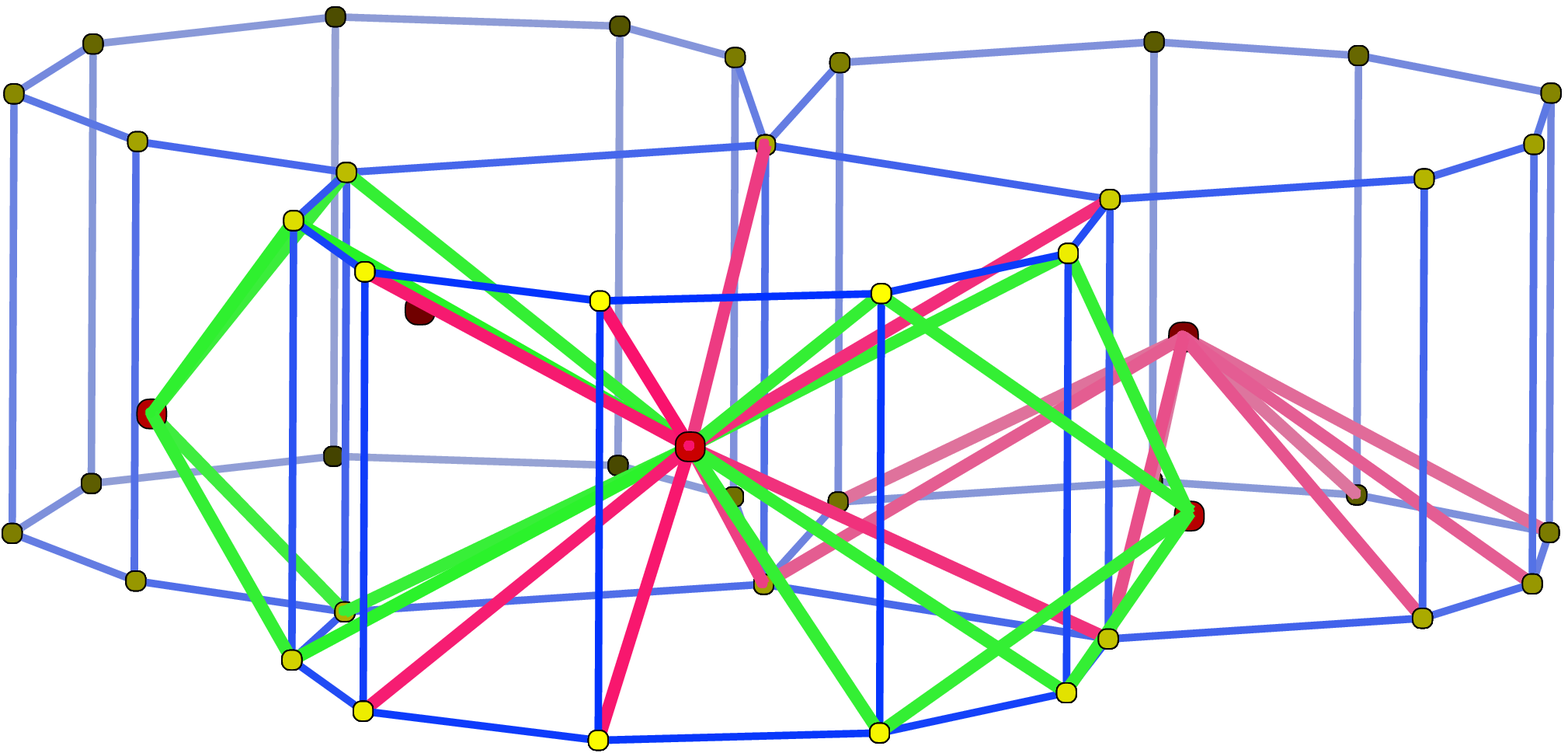}}
  \caption{The E-construction applied to $\Pi=\C_9^\alpha\times[0,1]$. \emph{Left:}
    Three of the nine prisms of~$\Pi$. \emph{Right:} Three of the 18
    octagonal pyramids and two out of ${9\choose 2}$ octahedra
    of~$E(\Pi)$.}
  \label{fig:eppstein}
\end{figure}

\begin{example}
  Let $\C$ be a cellulation of a closed orientable surface~$S$ with
  $f$-vector $(f_0,\,f_1,\,f_2)$. Then applying
  Proposition~\ref{prop:e-construction} to $\C\times[0,1]$ yields a
  cellulation of the prism $S\times[0,1]$ with $2f_0+f_2$ vertices
  consisting of $f_1$~octahedra and $2f_2$~pyramids; see
  Figure~\ref{fig:eppstein}.
\end{example}

\subsection{Heegaard splittings}

\begin{figure}[htbp]
  \centering
  \includegraphics[width=.5\linewidth]{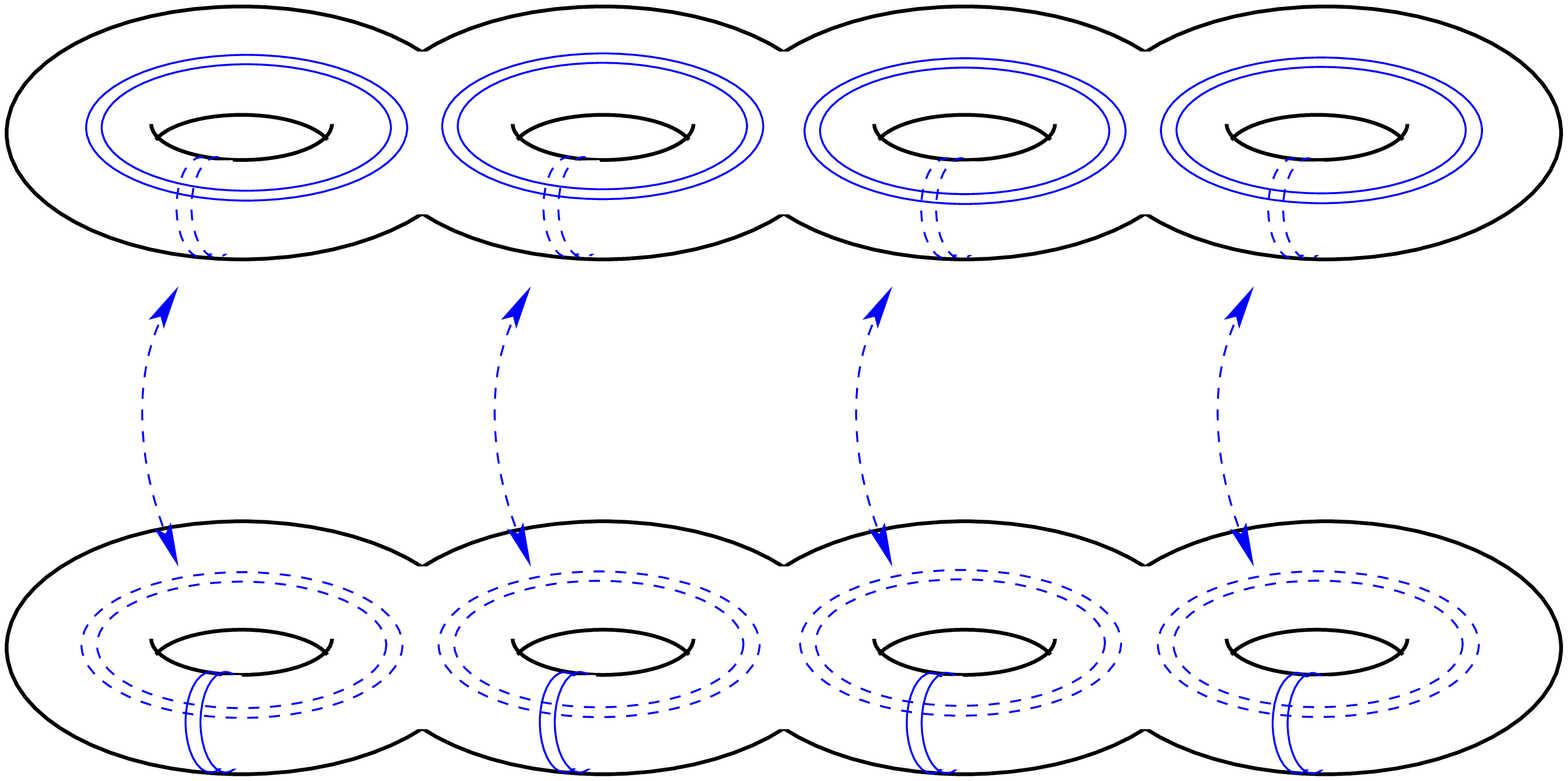}
  \caption{Heegaard splitting of $S^3$ of genus $g=4$. The complement of
    one solid handlebody in the $3$-sphere is the other solid
    handlebody of the same genus. One copy of each doubled solid
    (resp.~dashed) homology generator on the upper handlebody~$H_1$ is
    identified with one copy of the solid (resp.~dashed) one on the
    lower handlebody~$H_2$ in the way indicated by the arrows, and the
    union of all copies of the generators induces a cellulation
    of~$H_1$ (resp.~$H_2$) into one $3$-ball and $g$~solid cylinders.}
  \label{fig:heegaard}
\end{figure}

\begin{proposition}\label{prop:heegaard}
  \textup{(see \cite[Section 8.3.2]{Stillwell93})} For any $g\ge1$,
  the $3$-sphere may be decomposed into two solid handlebodies
  $S^3=H_1\cup H_2$ that are identified along a surface $S_g=H_1\cap
  H_2$ of genus~$g$. Conversely, any $3$-manifold can be split into
  handlebodies~$H_1$, $H_2$ and is determined up to homeomorphism by
  the images $h(m_1),\dots,h(m_{2g})$ on~$H_2$ of the canonical
  meridians $m_1,\dots,m_{2g}$ of~$H_1$ under the identification map
  $h:\partial H_1\to\partial H_2$.\hfill$\Box$
\end{proposition}

\begin{theorem}\label{thm:lazarus-etal}
  \textup{(Lazarus et al.~\cite[Theorem 1]{Lazarus-etal01})} Any
  triangulation $\T$ of a closed orientable surface of genus~$g$ with
  a total of $f=f_0+f_1+f_2$ cells can be  refined to a
  triangulation~$\tilde\T$ with $O(fg)$ vertices that contains
  representatives of the canonical homology generators in its
  $1$-skeleton. These representatives only intersect in a single
  vertex, and each one uses~$O(f)$ vertices and edges.  \hfill$\Box$
\end{theorem}

\begin{proof}[Idea of proof.]
Lazarus et al.\ present two algorithms that actually
\emph{compute} the canonical homology generators,
and from which the subdivision is easy to derive. Both algorithms
are ``optimal'' from a worst-case complexity point of view.

The first algorithm is inductive, removing one triangle at a
time from the surface in question and maintaining information
about the still unvisited part of the surface and
its collared boundary.

The second algorithm (based on Brahana \cite{Brahana1921})
starts with a maximal subgraph $G$
of the vertex-edge graph of the surface that has a
connected complement $\T{\setminus}G$, which is thus an open disk.
One derives generators for the fundamental group of $G$, which also
generate the fundamental group of $\T$.
These generators are then modified to yield canonical generators for
the fundamental group of $\T$.
\end{proof}

\section{Many triangulated 3-spheres}

\begin{figure}[htb]
  \centering
  \input{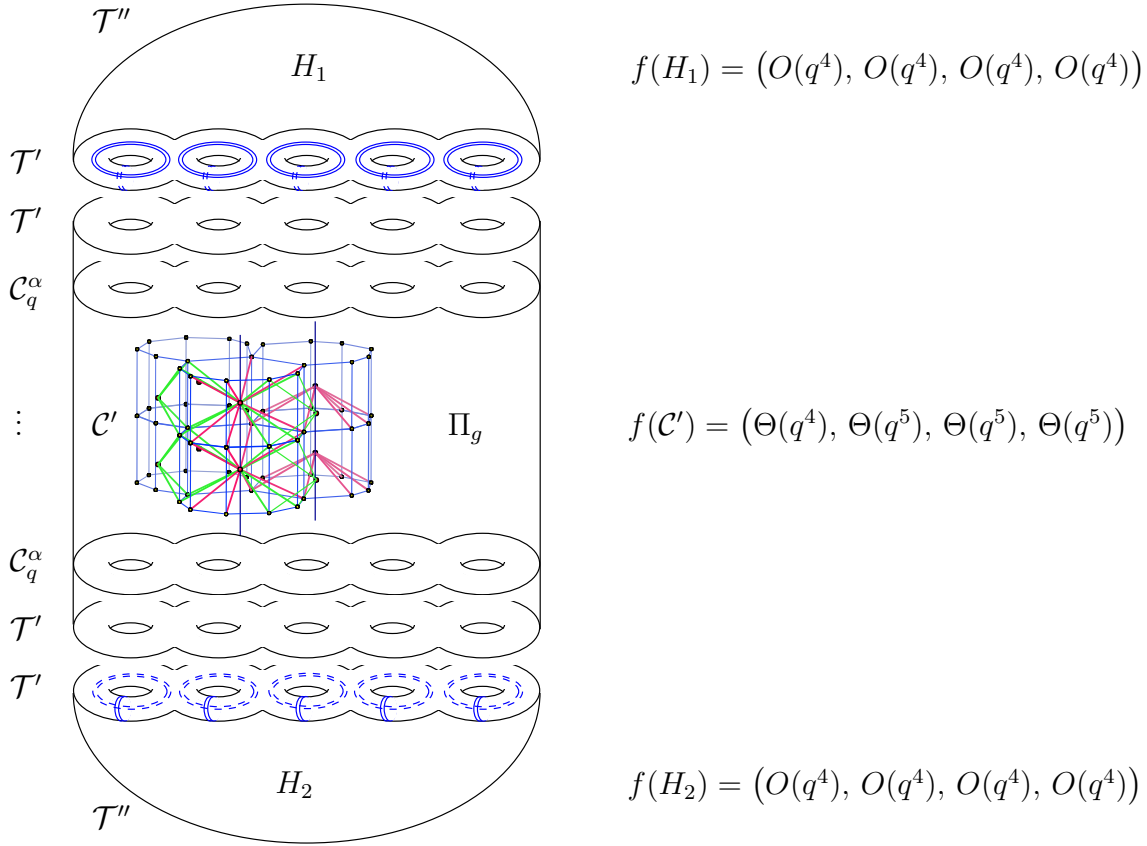}
  \caption{The thickened Heegaard splitting $S^3=H_1\cup\C'\cup H_2$
    of $S^3$. Not shown is the triangulation of the handlebodies~$H_1$
    and~$H_2$. Independently triangulating the $\Theta(n^{5/4})$
    octahedral $3$-cells of~$\C'$ in different ways yields ``many
    triangulated $3$-spheres''.}
  \label{fig:bigmac}
\end{figure}

\begin{proof}[Proof of Theorem~\ref{thm:many-spheres}:] We build a cellular
decomposition~$\S$ of~$S^3$ with $n$~vertices and $\Theta(n^{5/4})$
octahedral facets from two triangulated handlebodies and a stack of
prisms over a Heffter surface. The theorem then follows by
independently triangulating the octahedra.

The construction begins with a Heegaard splitting $S^3=H_1\cup H_2$
of~$S^3$ of genus $g=q(q-5)/4+1$ as in Proposition~\ref{prop:heegaard}, for
any prime power~$q$ of the form $q=4k+1$ for~$k\ge1$. We replace the
boundary $S_g=H_1\cap H_2$ of the handlebodies by the prism
$\Pi_g=S_g\times[0,1]$, pick a generator~$\alpha$ of~$\F_q^*$, and
embed a copy of the Heffter triangulation~$\T_q^\alpha$
on~$S_g\times\{0\}$ and~$S_g\times\{1\}$.
 
\begin{dense_itemize}
\item \emph{The triangulated handlebodies.} Use
  Theorem~\ref{thm:lazarus-etal} to refine each copy of~$\T_q^\alpha$
  to a triangulation of~$S_g$ that contains representatives of the
  canonical homology generators $\{a_i,b_i:1\le i\le g\}$ in its
  $1$-skeleton, such that the $a_i$'s span meridian disks in~$H_1$ and
  the $b_i$'s do the same in~$H_2$.  This introduces $O(q^2g)=O(q^4)$
  new vertices. Double all generators as in Figure~\ref{fig:heegaard}
  using another $O(q^4)$~vertices to obtain a triangulation~$\T'$
  of~$S_g$, and in each handlebody triangulate the meridian disks
  spanned by all these polygonal curves (using a total of
  $O(q^2g)$~triangles, but no new vertices). Then cone the boundary of
  each of the $2g$~solid cylinders bounded by the meridian disks to a
  new vertex (introducing a total of $O(q^2g)$ tetrahedra), and cone
  the triangulated boundary of each of the remaining two $3$-balls to
  another new vertex. This last step uses $2g+2$~new vertices and
  $O(q^4)$ tetrahedra.  The total $f$-vector of this
  triangulation~$\T''$ of~$H_1\cup H_2$ is
  \[
      \big(O(q^4),\,O(q^4),\,O(q^4),\,O(q^4)\big).
  \] 
   
\item \emph{The stack of prisms.} Let $\C_q^\alpha\times I_m$
  cellulate the manifold with boundary~$\Pi_g=S_g\times[0,1]$, where
  $I_m$~is the subdivision of~$[0,1]$ into $m=\Theta(q^3)$ closed
  intervals, and refine each of $\C_q^\alpha\times\{0\}$ and
  $\C_q^\alpha\times\{1\}$ into the triangulation~$\T'$. This refined
  cellulation~$\C$ is composed of $\Theta(q^4)$~prisms over
  $(q-1)$-gons and $2q$~\mbox{3-cells} whose boundary consists of
  $q-1$ $4$-gons, one $(q-1)$-gon, and on average $O(q^3)$ triangles
  that together triangulate another $(q-1)$-gon.  The boundary of~$\C$
  consists of the union of these $O(q^4)$ triangles, and its total
  $f$-vector is
  \[
     \big(\Theta(q^4),\,\Theta(q^5),\,\Theta(q^5),\,\Theta(q^4)\big).
  \]
  Apply the E-construction (Proposition~\ref{prop:e-construction})
  to~$\C$, using $\Theta(q^4)$~new vertices, to arrive at a
  cellulation~$E(\C)$ of~$\Pi_g$ into $\Theta(q)$~simplices (pyramids
  over the boundary triangles), $\Theta(q^4)$ bipyramids over
  $(q-1)$-gons, and $\Theta(q^5)$~octahedra. Now triangulate the
  bipyramids by joining each main diagonal to each edge of the base
  $(q-1)$-gon. This cellulation~$\C'$ of~$\Pi_g$ consists of
  $\Theta(q^5)$~simplices and $\Theta(q^5)$~octahedra
  (Figure~\ref{fig:bigmac}).  Its total $f$-vector is
  \[
      \big(\Theta(q^4),\,\Theta(q^5),\,\Theta(q^5),\,\Theta(q^5)\big).
  \]
\end{dense_itemize}

\noindent The desired cellulation of~$S^3$ is $\S=\T''\cup\C'$. 
\end{proof}

\medskip

\end{document}